\documentclass[12pt]{article}

\usepackage{morefloats}
    \usepackage{amsmath}
    \usepackage{amsfonts}
\usepackage{color}
\usepackage{multirow}
\usepackage{graphicx}
\usepackage{rotating}

    \newtheorem{theo}{Theorem}
    \newtheorem{coro}{Corollary}

    \def\0{{\bf 0}}

    \def\PP{\mathbb{P}}

    \def\epsilon{\varepsilon}

    \def\Ty0M{{\cal T}_{y_0}\M}

    \newcommand{\e}{{\textrm{e}}}
    \newcommand{\BE}{\mathbb{E}}
    \newcommand{\BV}{\mathbb{V}}
    \newcommand{\dd}{\textrm{d}}
    \newcommand{\QQ}{\mathbb{Q}}

\newcommand{\M}{{\cal M}}

    \newcommand{\Cov}{{\rm Cov}}

    \newcommand{\RR}{{\mathbb{R}}}
    \newcommand{\vertk}{{\stackrel{\cal D}{\longrightarrow}}}

    \def\bdm{\begin{displaymath}}
    \newcommand{\edm}{\end{displaymath}}
    \def\benu{\begin{enumerate}}
    \def\eenu{\end{enumerate}}
    \def\beqn{\begin{equation}}
    \def\eeqn{\end{equation}}
    \def\be{\begin{equation}}
    \def\ee{\end{equation}}
    \def\bea{\begin{eqnarray}}
    \def\eea{\end{eqnarray}}
    \newcommand{\bean}{\begin{eqnarray*}}
    \newcommand{\eean}{\end{eqnarray*}}
    \newcommand{\bear}{\begin{eqnarray}}
    \newcommand{\eear}{\end{eqnarray}}

    \def\qed{\hfill\hbox{${\vcenter{\vbox{
        \hrule height 0.4pt\hbox{\vrule width 0.4pt height 6pt
        \kern5pt\vrule width 0.4pt}\hrule height 0.4pt}}}$}}

\begin{document}

\title{\bf On a test of normality based on the empirical moment generating function}

\author{N. Henze and S. Koch}

\date{\today}
\maketitle

\footnotetext{ {\em American Mathematical Society 2000 subject
classifications.} Primary 62F05 Secondary 62G10} 

\begin{abstract} We provide the lacking theory for a test of normality based on the empirical moment generating function.
\end{abstract}

\noindent

\section{Introduction}\label{intro}
As evidenced by the recent papers of \cite{ghosh}, \cite{goia}, \cite{mbah}, \cite{szynal}, \cite{hanusz}, \cite{bouzeba}, \cite{desgagne}, \cite{shalit}, \cite{sanqui}, \cite{lee}, \cite{schick}, \cite{quessy}, \cite{harry} and
  \cite{goia2}, \cite{torabi}, there is an ongoing interest in testing of normality. This paper is not devoted to review the multitude of tests
suggested and studied for this testing problem (for an account of classical tests, see, e.g., \cite{thode} or \cite{dago}),  but to provide missing mathematical theory for a recent test suggested by \cite{zghoul}, which is based on the moment generating function.

To be specific, let $X_1,X_2, \ldots $ be independent and identically distributed (i.i.d.) random variables with an unknown continuous distribution, defined on a common probability space $(\Omega,{\cal A}, \mathbb{P})$.
Write $\mathbb{P}^{X_1}$ for the distribution of $X_1$, N$(\mu,\sigma^2)$ for the normal distribution with expectation $\mu$ and variance $\sigma^2$
 and ${\cal N} = \{{\rm{N}}(\mu,\sigma^2): \mu \in \mathbb{R}, \sigma^2 >0\}$ for the class of
(non-degenerate) normal distributions. Based on $X_1,\ldots,X_n$, \cite{zghoul} proposed to reject the hypothesis
\[
H_0: \mathbb{P}^{X_1} \in {\cal N}
\]
for large values of the test statistic
\begin{equation}\label{teststat}
T_{n,\beta}  = n \int_{-\infty}^\infty \left(M_n(t) - M_0(t) \right)^2 \, \exp\left(- \beta t^2\right) \, \textrm{d}t.
\end{equation}
Here, $\beta >2$ is a fixed parameter, and $M_0(t) = \exp(t^2/2)$, $t\in \mathbb{R}$, is the moment generating function of the standard normal distribution.
Moreover,
\[
M_n(t) = \frac{1}{n} \sum_{j=1}^n \exp\left(t Y_{n,j}\right), \quad t \in \mathbb{R},
\]
is the empirical moment generating function of the scaled residuals
\[
Y_{n,j} = \frac{X_j - \overline{X}_n}{S_n}, \quad j =1,\ldots,n,
\]
where $\overline{X}_n = n^{-1}\sum_{j=1}^n X_j$ stands for  the sample mean, and
\[
S_n^2 = \frac{1}{n} \sum_{j=1}^n \left(X_j - \overline{X}_n\right)^2
\]
denotes the sample variance of $X_1,\ldots,X_n$.

The rationale for considering $T_{n.\beta}$ as a genuine test statistic for normality is clear-cut: Under $H_0$, the standardized residuals $Y_{n,1},\ldots,Y_{n,n}$ should be, at least for large $n$,
approximately standard normally distributed. Hence, $M_n$ should be close to $M_0$, and some measure of deviation between $M_n$ and $M_0$ should yield a reasonable test statistic.
 Notice that $T_{n,\beta}$ is a weighted $L^2$-type statistic. Such statistics have been employed in
numerous goodness-of-fit testing problems (see, e.g. \cite{bahe}). If, in (\ref{teststat}),
  one replaces $M_n$ by the empirical {\em characteristic} function of $X_1,\ldots,X_n$ and $M_0(t)$ by $\exp(-t^2/2)$, the
characteristic function of the standard normal distribution, one obtains the statistic of \cite{epps}. For goodness-of-fit tests based on the empirical moment generating function, see, e.g.,
\cite{cabana}, \cite{eppssi}, \cite{kallioras}, \cite{meintanis}, and \cite{meintanis2}.

\vskip.3cm

Straightforward computation of the integral figuring in (\ref{teststat}) shows that $T_{n,\beta}$ takes the form
\begin{equation}\label{sumdarst}
T_{n,\beta}  =  \sqrt{\pi}\bigg(\! \frac{n}{\sqrt{\beta \! -\! 1}} -\frac{2}{\sqrt{\beta\! -\! \frac{1}{2}}} \sum_{i=1}^{n} \! \exp\! \Big(\frac{Y_{n,i}^2}{4\beta \! -\! 2}\Big)
+ \frac{1}{n\sqrt{\beta}} \! \sum_{i,j=1}^{n} \! \exp \! \Big(\frac{(Y_{n,i}\! +\! Y_{n,j})^2}{4\beta}\Big)\bigg),
\end{equation}
which is amenable to computational purposes. A simulation study conducted by \cite{zghoul} showed that the test based on $T_{n,\beta}$ is a strong competitor to classical tests of normality, such as the
Anderson-Darling test, the Shapiro-Wilk test, the Epps-Pulley test, and the D'Agostino test (for an account of these procedures, see \cite{dago}).

The purpose of this paper is to provide some theoretical background for the test of Zghoul. We will prove that $T_{n,\beta}$ has a non-degenerate limit distribution under $H_0$, and we will show that the test
is consistent against general alternatives. Moreover, letting the parameter $\beta$ tend to infinity, $T_{n,\beta}$ approaches, upon suitable centering and rescaling,
squared sample skewness, which is one of the first statistics used for testing for normality.

 The rest of the paper is organized as follows. In Section \ref{secnull}, we state a result on the limit null distribution of $T_{n,\beta}$ and derive the expectation and the variance of this
 limit law. Section \ref{consist}  is devoted to the behavior of $T_{n,\beta}$
under a fixed alternative to normality, and Section \ref{betainfti} considers the case $\beta \to \infty$. Some technical proofs are deferred to Section \ref{secasymneg}. The paper concludes with some
remarks and open problems.

\section{The limit null distribution of $T_{n,\beta}$}\label{secnull}
In this section, we derive the limit distribution of $T_{n,\beta}$ under $H_0$. Since $T_{n,\beta}$ is invariant with respect to affine transformations of $X_1,\ldots,X_n$, the null distribution
of $T_{n,\beta}$ does not depend on the true values of $\mu$ and $\sigma^2$. We thus assume without loss of generality that $\mu=0$ and $\sigma^2=1$ throughout this section.
Since $T_{n,\beta}$ is a weighted $L^2$-statistic, a convenient setting for asymptotics is the separable Hilbert space
${\cal H} = L^2(\RR,{\cal B},w(t)\textrm{d}t)$ of (equivalence classes of) measurable functions $f:\RR \to \RR$ such that $\int_\RR f^2(t) w(t)\, \textrm{d} t < \infty$.
Here, ${\cal B}$ is  the $\sigma$-field of Borel sets of $\RR$, and $w(t) = \exp(-\beta t^2)$. The inner product and the resulting norm on ${\cal H}$ will be denoted by
\[
\langle f,g \rangle = \int_{\RR} f(t)g(t) \, w(t) \, \textrm{d}t, \quad \|f \| = \left(\int_{\RR} f^2(t) \, w(t)\, \textrm{d}t \right)^{1/2},
\]
respectively.  Putting
\begin{equation}\label{defwnt}
W_n(t) = \sqrt{n} \left(M_n(t) - M(t) \right), \quad t \in \RR,
\end{equation}
$W_n$ is a random element of ${\cal H}$, and we have $T_{n,\beta} = \|W_n\|^2$. If we could prove
$W_n \vertk W$ in ${\cal H}$ for some random element $W$ of ${\cal H}$, where $\vertk$ denotes convergence in distribution in ${\cal H}$,
the continuous mapping theorem would yield $T_{n,\beta} \vertk \|W\|^2$. If not stated otherwise, convergence is always meant as $n \to \infty$.

\medskip

\begin{theo}\label{nulldistri}
There is a centred Gaussian element $W$ of ${\cal H}$ having covariance kernel
\begin{equation}\label{kernel}
K(s,t) = \e^{(t^2+s^2)/2}\bigg(\e^{ts} - 1 -ts-\frac{t^2s^2}{2}\bigg), \quad s,t \in \RR,
\end{equation}
such that $W_n \, \vertk \, W$.
\end{theo}

\medskip

\begin{coro}\label{cor1} Under $H_0$, we have
\[
T_{n,\beta} \, \vertk \, \|W\|^2 = \int_{\RR} W^2(t) \, \e^{-\beta t^2} \, {\rm d} t,
\]
where $W$ is the Gaussian element of ${\cal H}$ figuring in Theorem \ref{nulldistri}.

\end{coro}

\medskip

{\sc Proof.} The main problem in proving Theorem \ref{nulldistri} is that $nM_n(t) = \sum_{j=1}^n \exp(tY_{n,j})$ is not a sum of i.i.d.
random variables. To overcome this drawback, notice that
\[
\e^{tY_{n,i}} - \e^{tX_i} = \e^{X_i}\big(\e^{t(Y_{n,i} -X_i)}-1\big),
\]
where
\[
Y_{n,i}-X_i = \frac{X_i(1-S_n)-\overline{X}_n}{S_n}.
\]
Taylor's theorem yields
\[
\e^{t(Y_{n,i} -X_i)} - 1 = t(Y_{n,i}-X_i)+ \frac{1}{2} t^2(Y_{n,i} -X_i)^2 \exp\big(\Theta_{n,i}t(Y_{n,i}-X_i)\big).
\]
Here, $\Theta_{n,i} = \Theta_{n,i}(t,X_i,X_1, \dots, X_n)$ are random variables with $|\Theta_{n,i}|\leq1$.

It follows that
\begin{align*}
\frac{1}{\sqrt{n}} \sum_{i=1}^{n}\e^{tY_{n,i}} - \frac{1}{\sqrt{n}} \sum_{i=1}^{n}\e^{tX_i}
= \frac{1}{\sqrt{n}} \sum_{i=1}^{n} \e^{tX_i}\, t \, \frac{X_i(1-S_n)-\overline{X}_n}{S_n}  +R_n(t),
\end{align*}
where
\[
R_n(t) :=  \frac{1}{\sqrt{n}} \sum_{i=1}^{n}\e^{tX_i} \, \frac{t^2}{2}\left(\!  \frac{X_i(1\! -\! S_n)\! -\! \overline{X}_n}{S_n}\!  \right)^2\!
 \exp\bigg(\! \Theta_{n,i}t\frac{X_i(1\! -\! S_n)\! -\! \overline{X}_n}{S_n}\bigg).
\]
The main part of the proof consists of showing
\begin{equation}\label{asymneg}
\|R_n\|^2 \, = \ o_\PP(1).
\end{equation}
Since the proof of (\ref{asymneg}) is quite technical due to the unboundedness of the moment generating funktion
over the whole line, it is deferred to Section \ref{secasymneg}.

Since $S_n = 1 + o_\PP(1)$ (remember that $\mu=0$ and $\sigma^2 =1$), we have
\begin{align*}
\frac{1}{\sqrt{n}}\sum_{i =1}^{n}\e^{tX_i}t  &\frac{X_i(1-S_n)-\overline{X}_n}{S_n} \\
&= \frac{(1-S_n^2)}{(1+S_n)S_n}\cdot \frac{1}{\sqrt{n}}\sum_{i =1}^{n}tX_i\e^{tX_i} - \frac{\overline{X}_n}{S_n}\cdot \frac{1}{\sqrt{n}}\sum_{i =1}^{n}t\e^{tX_i}\\
&= \frac{(1-S_n^2)}{2}\cdot \frac{1}{\sqrt{n}}\sum_{i =1}^{n}tX_i\e^{tX_i} - \overline{X}_n\cdot \frac{1}{\sqrt{n}}\sum_{i =1}^{n}t\e^{tX_i} + r_{n,1}(t),
\end{align*}
where $r_{n,1}$ is a random element of ${\cal H}$ satisfying $\|r_{n,1}\| = o_\PP(1)$. Now, use
\[
\sqrt{n} \left(S_n^2-1 \right) = \frac{1}{\sqrt{n}} \sum_{j=1}^n \left(X_j^2 - 1 \right) + o_\PP(1)
\]
to show that
\[
\frac{(1-S_n^2)}{2}\cdot \frac{1}{\sqrt{n}}\sum_{i =1}^{n}tX_i\e^{tX_i} = -\frac{1}{2\sqrt{n}}\cdot\frac{1}{n} \sum_{i,j = 1}^{n}(X_j^2-1)tX_i\e^{tX_i} + r_{n,2}(t),
\]
where the random element $r_{n,2}$ of ${\cal H}$ satisfies $\|r_{n,2}\| = o_\PP(1)$.
Next, let
\[
E_1(t) := \BE\left[ tX_1\e^{tX_1} \right] = t^2\e^{t^2/2}, \quad
E_2(t) := \BE\left[ t\e^{tX_1} \right] = t\e^{t^2/2}, \quad t \in \RR,
\]
and invoke the law of large numbers to end up in
\[
\frac{1}{\sqrt{n}}\sum_{i =1}^{n}\e^{tX_i}t  \frac{X_i(1-S_n)-\overline{X}_n}{S_n} = -\frac{1}{\sqrt{n}}\sum_{j = 1}^{n}\bigg((X_j^2-1)\frac{E_1(t)}{2} + X_jE_2(t) \bigg)+r_{n,3}(t),
\]
where $r_{n,3} \in  {\cal H}$ and $\|r_{n,3}\| = o_\PP(1)$.
Putting
\begin{equation}\label{hvonxt}
h(x,t)= \e^{tx} -\e^{t^2/2} -(x^2-1)\frac{E_1(t)}{2} - xE_2(t), \quad x,t \in \RR,
\end{equation}
and
\begin{equation}\label{wntilde}
\widetilde{W}_n(t) :=  \frac{1}{\sqrt{n}}\sum_{i=1}^{n}h(X_i,t),
\end{equation}
the definition of $W_n$ (see (\ref{defwnt})), the reasoning given above and (\ref{asymneg}) imply
\begin{equation}\label{wnwntild}
W_n(t) = \widetilde{W}_n(t) + \Delta_n(t),
\end{equation}
where $\Delta_n$ is a random element of ${\cal H}$ satisfying $\|\Delta_n\| = o_\PP(1)$.

Now, some algebra yields $\BE h(X_1,t) = 0$, $t \in \RR$, and $\BE[h(X_1,s)h(X_1,t)] = K(s,t)$, $s,t \in \RR$,
where $K$ is given in (\ref{kernel}). Since the random elements $h(X_j, \cdot)$, $j=1,\ldots,n$, of ${\cal H}$ figuring in
(\ref{wntilde}) are i.i.d., a Hilbert space central limit theorem (see. eg., Theorem 1.1. of \cite{kundu}) gives
$\widetilde{W}_n \vertk W$, where $W$ is a centred Gaussian element of ${\cal H}$ having covariance kernel $K$.
In view of $\|W_n - \widetilde{W}_n\| = o_\PP(1)$, it follows that $W_n \vertk W$.

\vskip.5cm

It is well-known that the distribution of
\[ T_\infty  := \|W\|^2
\]
is that of $\sum_{j \ge 1} \lambda_j N_j^2$, where $N_1,N_2, \ldots $ are i.i.d. standard normal random variables,
and $\lambda_1,\lambda_2, \ldots $ are the nonzero eigenvalues corresponding to the orthonormal eigenfunctions
of the integral operator $A: {\cal H} \rightarrow {\cal H}$, where
\[
(Af)(t) = \int_\RR K(s,t) f(s) \, \exp(-\beta s^2) \, \dd s , \quad f \in {\cal H},
\]
and $K$ is given in (\ref{kernel}).
We did not succeed in solving this integral equation. However, using formulae of \cite{shorack}, p. 213, we obtain the following information on the distribution of $T_\infty$.

\medskip

\begin{theo} We have
\begin{enumerate}
\item[a)] $\BE(T_\infty) = \displaystyle{\frac{\sqrt{\pi}}{\sqrt{\beta -2}} - \frac{\sqrt{\pi}}{\sqrt{\beta -1}}\left(1+\frac{1}{2(\beta -1)} + \frac{3}{8(\beta -1)^2}\right)}$,\\
\item[b)]  $\BV\left(T_\infty\right) = 2\pi \! \displaystyle{\bigg( \! \frac{1}{\sqrt{\beta}\sqrt{\beta \! - \! 2}} - \frac{4}{\sqrt{\gamma}}-\frac{6}{\gamma^{3/2}}-\frac{6}{\gamma^{5/2}} +  \frac{1}{\beta \! -\! 1} + \frac{1}{2(\beta \! - \! 1)^3} + \frac{9}{64(\beta \! - \! 1)^5}\! \bigg)}$,\\
    where $\gamma = 4(\beta -1)^2-1$.
\end{enumerate}
\end{theo}

\medskip

{\sc Proof:} Since
\[
\BE(T_\infty) = \int_\RR K(t,t) w(t)\,  \dd t
\]
and
\[ \BV(T_\infty) = 2 \int_{\RR^2} K^2(s,t) w(s)w(t)\,  \dd s \dd t,
\]
the result follows from tedious but straightforward calculations of integrals.

\section{Consistency}\label{consist}
In this section we show that the test for normality based on $T_{n,\beta}$ is consistent against general alternatives.
Our main result is as follows.

\medskip

\begin{theo}\label{fixedalt}
Assume that $X_1$ has a non-degenerate distribution, and that the moment generating function $M(t) = \BE \exp(tX_1)$ exists for each $t \in \RR$.
We then have
\begin{equation}\label{unglie}
\liminf_{n \to \infty} \frac{T_{n,\beta}}{n} \ge \int_\RR \left( M(t) - M_0(t)\right)^2 \e^{-\beta t^2}\, \dd t \quad \PP\textrm{-a.s.}.
\end{equation}
\end{theo}

\medskip

{\sc Proof.}
Remember that $X_1,X_2, \ldots $ are defined on the probability space $(\Omega,{\cal A},\PP)$. In view of affine-invariance, we assume w.l.o.g. $\BE(X_1) =0$ and
$\BV(X_1) =1$. Fix $\varepsilon \in (0,1)$.
By the Strong Law of Large Numbers there is a set $\Omega_0 =  \Omega_0(\varepsilon) \in \mathcal{A}$ with $\PP(\Omega_0) = 1$ such that, for each $\omega \in \Omega_0$
there is an integer $n_0 = n_0(\varepsilon)$ and
\[ |\overline{X}_n(\omega)| \leq \varepsilon, \  |S_n(\omega)-1| \leq \varepsilon
\]
for each $n \ge n_0$. Putting $M_n(t,\omega):= n^{-1}\sum_{i=1}^{n}\exp(tX_i(\omega)),\ \omega \in \Omega, \ t \in \RR$, we obtain for each $n\geq n_0(\omega)$ and $t \geq 0$,
\begin{align*}
\frac{1}{n}\sum_{i = 1}^{n}\exp\bigg(t\frac{X_i(\omega) -\varepsilon}{1+\varepsilon}\bigg) \leq M_n(t,\omega) \leq \frac{1}{n}\sum_{i = 1}^{n}\exp\bigg(t\frac{X_i(\omega) +\varepsilon}{1-\varepsilon}\bigg).
\end{align*}
Again by the Strong Law of Large Numbers there is a set $\Omega_1 = \Omega_1(\varepsilon,t)$ depending on $\varepsilon$ and $t$
 with $\PP(\Omega_1)=1$ such that for each $\omega \in \Omega_1:$
\begin{align*}
\BE\left[\exp\bigg(t\frac{X_1-\varepsilon}{1+\varepsilon}\bigg)\right] \leq \varliminf _{n \to \infty}M_n(t,\omega) \leq \varlimsup _{n \to \infty}M_n(t,\omega) \leq \BE\left[\exp\bigg(t\frac{X_1 +\varepsilon}{1-\varepsilon}\bigg)\right].
\end{align*}
Letting $\varepsilon \downarrow 0$ then yields $M_n(t,\cdot) \to M(t)$ $\PP$-almost surely for fixed $t\geq 0$.

If $t <0$, we have for each $n\geq n_0(\omega)$
\begin{align*}
\frac{1}{n}\sum_{i = 1}^{n}\exp\bigg(t\frac{X_i(\omega) +\varepsilon}{1-\varepsilon}\bigg) \leq M_n(t,\omega) \leq \frac{1}{n}\sum_{i = 1}^{n}\exp\bigg(t\frac{X_i(\omega) -\varepsilon}{1+\varepsilon}\bigg),
\end{align*}
and the same reasoning entails $M_n(t) \to M(t)$ almost surely for each fixed $t\in \RR$.
In other words, for each $t \in \RR$ there is a set $\Omega_2(t) \in \mathcal{A}$ with $\PP(\Omega_2(t))=1$ and
\[
M_n(t,\omega) \rightarrow M(t) \text{ for each } \omega \in \Omega_2(t).
\]
Writing $\QQ$ for the set of rational numbers, it follows that
\[
M_n(t,\omega) \rightarrow M(t) \ \forall t \in \QQ
\]
for each $\omega \in \Omega_3 := \bigcap_{t \in \QQ} \Omega_2(t)$.
Since $M_n$ and $M$ are convex functions and $\QQ$ is dense in $\RR$, we have for each $\omega \in \Omega_3$ that $M_n(t,\omega) \rightarrow M(t), \ t \in I,$ where $I$ is an arbitrary compact set and thus
\begin{align*}
M_n(t,\omega) \rightarrow M(t), \ t \in \RR,
\end{align*}
for each $\omega \in \Omega_3$ (e.g. see \cite{roberts} C.7, p. 20). Now fix $\omega \in \Omega_3$.
By Fatou's lemma,
\begin{eqnarray*}
\varliminf_{n \to \infty} \frac{T_{n,\beta}(\omega)}{n} & = & \varliminf_{n \to \infty} \int_\RR \left(M_n(t,\omega)-M_0(t)\right)^2 \e^{-\beta t^2} \, \dd t\\
 &\geq &  \int_\RR \varliminf_{n \to \infty} \left(M_n(t,\omega)-M_0(t)\right)^2 \e^{-\beta t^2} \, \dd t\\
& = & \int_\RR \left(M(t)-M_0(t)\right)^2 \e^{-\beta t^2} {\rm d} t,
\end{eqnarray*}
as was to be shown.

\vskip.5cm

If the distribution of $X_1$ is non-normal and satisfies the conditions of Theorem \ref{fixedalt}, the right-hand side of (\ref{unglie}) is strictly positive,
and thus $T_n \to \infty$ $\PP$-a.s. Therefore, due to Corollary \ref{cor1}, the test for normality based on $T_{n,\beta}$ is consistent against any such alternative.

\section{The case $\beta \to \infty$}\label{betainfti}
In this section we analyse the asymptotic behaviour of the test statistic $T_{n,\beta}$ for fixed $n$ and $\beta \to \infty$. It will be seen that,
after a suitable centering and scaling, $T_{n,\beta}$ approaches the square of the first nonzero component of Neyman's smooth test for normality, which is squared sample skewness.
For an account on smooth tests of fit, see \cite{rayner}.

\begin{theo}\label{betainfty}
We have
\[
\lim_{\beta \to \infty} \frac{96}{5}\beta^{7/2} \left(\frac{T_{n,\beta}}{n\sqrt{\pi}} - \tau(\beta)\right) = b_{n,1}^2,
\]
where
\[
\tau(\beta) = \frac{1}{\sqrt{\beta-1}} -\frac{2}{\sqrt{\beta-\frac{1}{2}}} -\frac{2}{(4\beta-2)\sqrt{\beta- \frac{1}{2}}}  +\frac{1}{\sqrt{\beta}} + \frac{1}{2\beta^{3/2}}+\frac{3}{16\beta^{5/2}}
\]
and
\begin{equation}\label{skewness}
b_{n,1}  =  \frac{\frac{1}{n}\sum_{i = 1}^{n}(X_i-\overline{X}_n)^3}{S_n^3}
\end{equation}
denotes sample skewness of $X_1,\ldots,X_n$.
\end{theo}

{\sc Proof.}
We start with (\ref{sumdarst}) and notice that the scaled residuals $Y_{n,i}$ satisfy
\[
\sum_{i =1}^n Y_{n,i}    =  0, \quad \sum_{i =1}^n Y_{n,i}^2 = 1, \quad
\sum_{i =1}^n Y_{n,i}^3  =   nb_{n,1},\quad \sum_{i =1}^n Y_{n,i}^4 = nb_{n,2},
\]
where $b_{n,1}$ is given in (\ref{skewness}) and
\[
b_{n,2}  =  \frac{\frac{1}{n}\sum_{i = 1}^{n}(X_i-\overline{X}_n)^4}{S_n^4}
\]
is sample kurtosis of $X_1,\ldots,X_n$. Expanding the exponential terms figuring in
(\ref{sumdarst}) we have
\begin{align*}
\sum_{i = 1}^{n} \exp\bigg(\frac{Y_{n,i}^2}{4\beta-2}\bigg) &= \sum_{i = 1}^{n} \bigg(1+ \frac{Y_{n,i}^2}{4\beta-2} +  \frac{Y_{n,i}^4}{2(4\beta-2)^2} + \frac{Y_{n,i}^6}{6(4\beta-2)^3} + O\big(\beta^{-4}\big) \bigg)\\
& =n + \frac{n}{4\beta-2} +\frac{1}{2(4\beta-2)^2}nb_{n,2} + \frac{1}{6(4\beta-2)^3} \sum_{i=1}^{n}Y_{n,i}^6 + O\big(\beta^{-4}\big)
\end{align*}
and
\begin{align*}
\sum_{i,j = 1}^{n} \exp\bigg(\frac{(Y_{n,i}+Y_{n,j})^2}{4\beta}\bigg) &= \sum_{i,j = 1}^{n} \bigg(1+ \frac{(Y_{n,i}+Y_{n,j})^2}{4\beta} +  \frac{(Y_{n,i}+Y_{n,j})^4}{32\beta^2} + \frac{(Y_{n,i}+Y_{n,j})^6}{384\beta^3}\\
& \hspace{2cm} + O\big(\beta^{-4}\big)\bigg)\\
&= n^2 + \frac{n^2}{2\beta}+\frac{n^2}{16\beta^2}b_{n,2}+\frac{3n^2}{16\beta^2}+\frac{n}{192\beta^3}\sum_{i=1}^{n}Y_{n,i}^6 +\frac{5n^2}{64\beta^3}b_{n,2}\\
& + \frac{5n^2}{96\beta^3}b_{n,1}^2 + O\big(\beta^{-4}\big).
\end{align*}
Since
\[
\frac{1}{6(4\beta -2)^3} = \frac{1}{384\beta^3} + O\left(\beta^{-4}\right),
\]
it follows that
\begin{align*}
\frac{T_{n,\beta}}{\sqrt{\pi}}-\frac{n}{\sqrt{\beta-1}}& = -\frac{2}{\sqrt{\beta-\frac{1}{2}}}\bigg(n + \frac{n}{4\beta-2} +\frac{1}{2(4\beta-2)^2}nb_{n,2} + \frac{1}{384\beta^3}\sum_{i=1}^{n}Y_{n,i}^6\bigg)\\
&+\frac{1}{n\sqrt{\beta}}\bigg( n^2 + \frac{n^2}{2\beta}+\frac{n^2}{16\beta^2}b_{n,2}+\frac{3n^2}{16\beta^2}+\frac{n}{192\beta^3}\sum_{i=1}^{n}Y_{n,i}^6 +\frac{5n^2}{64\beta^3}b_{n,2}\\
& + \frac{5n^2}{96\beta^3}b_{n,1}^2 \bigg)+O\big(\beta^{-9/2}\big)
\end{align*}
and hence
\begin{align*}
\frac{T_{n,\beta}}{\sqrt{\pi}}& -\frac{n}{\sqrt{\beta-1}}  +\frac{2n}{\sqrt{\beta-\frac{1}{2}}} +\frac{2n}{(4\beta-2)\sqrt{\beta- \frac{1}{2}}} -\frac{n}{\sqrt{\beta}} - \frac{n}{2\beta^{3/2}}-\frac{3n}{16\beta^{5/2}}\\
& = \bigg(\frac{1}{192\beta^{7/2}} - \frac{1}{192\beta^3 \sqrt{\beta-\frac{1}{2}}} \bigg)\sum_{i=1}^{n}Y_{n,i}^6 +\frac{5n}{96\beta^{7/2}}b_{n,1}^2\\
& + \bigg(\frac{n}{16\beta^{5/2}} -\frac{n}{(4\beta-2)^2\sqrt{\beta-\frac{1}{2}}}\bigg)b_{n,2} + \frac{5n}{64\beta^{7/2}}b_{n,2} + o\big(\beta^{-7/2}\big).
\end{align*}
Since
\[
\frac{1}{192\beta^3 \sqrt{\beta-\frac{1}{2}}} = \frac{1}{192\beta^{7/2}} + o\left(\beta^{-7/2}\right)
\]
and
\begin{align*}
\frac{n}{16\beta^{5/2}} -\frac{n}{(4\beta-2)^2\sqrt{\beta-\frac{1}{2}}} &= \frac{n}{16\beta^{5/2}} - \frac{n}{16\beta^2\left(1-\frac{1}{2\beta}\right)^2\sqrt{\beta}\sqrt{1-\frac{1}{2\beta}}}\\
&= -\frac{5n}{64\beta^{7/2}} + o\left(\beta^{-7/2}\right),
\end{align*}
the result follows from
\begin{align*}
\frac{T_{n,\beta}}{\sqrt{\pi}}&-\frac{n}{\sqrt{\beta-1}} +\frac{2n}{\sqrt{\beta-\frac{1}{2}}} +\frac{2n}{(4\beta-2)\sqrt{\beta- \frac{1}{2}}}  -\frac{n}{\sqrt{\beta}} - \frac{n}{2\beta^{3/2}}-\frac{3n}{16\beta^{5/2}}\\
& =\frac{5n}{96\beta^{7/2}}b_{n,1}^2 + o\left(\beta^{-7/2}\right).
\end{align*}

Notice that Theorem \ref{betainfty} corresponds to Theorem 3.1 of \cite{baguhe} for the Epps-Pulley test statistic.\\

\section{The proof of (\ref{asymneg})}\label{secasymneg}

Since  $|\Theta_{n,i}|\leq  1$ and $(a \pm b)^2 \leq 2a^2 +2b^2$, for $a,b\in \RR$, we have
\begin{align*}
0 \leq R_n(t) \leq R_{n,1}(t) + R_{n,2}(t),
\end{align*}
where
\begin{align*}
R_{n,1}(t) &= \frac{(1-S_n)^2}{S_n^2}\cdot \frac{1}{\sqrt{n}} \sum_{i=1}^{n}\e^{tX_i}t^2X_i^2 \exp\bigg(|t|\frac{|X_i(1-S_n)-\overline{X}_n|}{S_n}\bigg),\\
R_{n,2}(t) &= \frac{\overline{X}_n^2}{S_n^2}\cdot \frac{1}{\sqrt{n}} \sum_{i=1}^{n}\e^{tX_i} t^2 \exp\bigg(|t|\frac{|X_i(1-S_n)-\overline{X}_n|}{S_n}\bigg).
\end{align*}
This decomposition yields $R_n^2(t) \leq 2R_{n,1}^2(t) + 2R_{n,2}^2(t)$ and thus
$$\|R_n\|^2 \leq 2\|R_{n,1}\|^2 + 2\|R_{n,2}\|^2.$$
Since
$$\left(\frac{(1-S_n)^2}{\sqrt{n}S_n^2}\right)^2 = O_\PP\left(n^{-3}\right), \qquad  \left(\frac{\overline{X}_n^2}{\sqrt{n}S_n^2}\right)^2 = O_\PP\left(n^{-3}\right),$$
we have
\begin{align*}
\|R_{n,1}\|^2 &\leq O_\PP\left(n^{-3}\right)\\
& \cdot\sum_{i,j=1}^{n}X_i^2X_j^2 \int_\RR \e^{t(X_i+X_j)} t^4 \exp\bigg( \frac{|t|}{S_n}\big((|X_i|+|X_j|)|1-S_n|+2|\overline{X}_n| \big) \bigg)  \e^{-\beta t^2}\dd t,\\
\|R_{n,2}\|^2 & \leq O_\PP\left(n^{-3}\right)\\
&  \cdot\sum_{i,j=1}^{n}\int_\RR \e^{t(X_i+X_j)} t^4 \exp\bigg( \frac{|t|}{S_n}\big((|X_i|+|X_j|)|1-S_n|+2|\overline{X}_n| \big) \bigg)  \e^{-\beta t^2}\dd t.
\end{align*}
Putting
\[ \alpha_n: = \alpha_n(i,j) := \frac{(|X_i|+|X_j|)|1-S_n|+2|\overline{X}_n|}{S_n}
\]
and observing that
\begin{eqnarray*}
(X_i+X_j \pm \alpha_n)^2 & \leq &  2(X_i+X_j)^2 + 2\alpha_n^2 ,\\
(X_i+X_j \pm \alpha_n)^4 & \leq &  4(X_i+X_j)^4 + 4\alpha_n^4,
\end{eqnarray*}

we obtain
\begin{align*}
&\int_\RR t^4 \exp\big(-\beta t^2 + t(X_i+X_j) + \alpha_n |t| \big) \dd t\\
& = \int_{0}^{\infty} t^4 \exp\big(-\beta t^2 + t(X_i+X_j) + \alpha_n t \big) \dd t\\
& + \int_{-\infty}^{0} t^4 \exp\big(-\beta t^2 + t(X_i+X_j) - \alpha_n t \big) \dd t\\
& \leq  \int_\RR t^4 \exp\big(-\beta t^2 + t(X_i+X_j + \alpha_n) \big) \dd t\\
& + \int_\RR t^4 \exp\big(-\beta t^2 + t(X_i+X_j - \alpha_n) \big) \dd t\\
& = \frac{\sqrt{\pi}\big((X_i+X_j + \alpha_n)^4+12\beta(X_i+X_j + \alpha_n)^2 + 12\beta^2\big)}{16\beta^{9/2}}\exp\bigg(\frac{(X_i+X_j + \alpha_n)^2}{4\beta}\bigg)\\
&+  \frac{\sqrt{\pi}\big((X_i+X_j - \alpha_n)^4+12\beta(X_i+X_j - \alpha_n)^2 + 12\beta^2\big)}{16\beta^{9/2}}\exp\bigg(\frac{(X_i+X_j - \alpha_n)^2}{4\beta}\bigg)\\
&\leq \frac{\sqrt{\pi}}{4\beta^{9/2}}\big( (X_i+X_j)^4 + \alpha_n^4+ 6\beta(X_i+X_j)^2 + 6\alpha_n^2 + 3\beta^2\big)\\
&\hspace{3cm}\cdot \left[\exp\bigg(\frac{(X_i+X_j + \alpha_n)^2}{4\beta}\bigg)+\exp\bigg(\frac{(X_i+X_j - \alpha_n)^2}{4\beta}\bigg)\right]
.
\end{align*}
Defining
\begin{eqnarray*}
C_n & := & \frac{2\max_{1\leq i\leq n} \{ |X_i| \} |1-S_n|+2|\overline{X}_n|}{S_n},\\
D_n & := & 2\max_{1\leq i\leq n} \{ |X_i| \}\cdot C_n,
\end{eqnarray*}
it follows that $\alpha_n(i,j) \le C_n$ and $|(X_i+X_j)\alpha_n(i,j)| \le D_n$ and thus
\begin{align*}
\begin{split}
\exp\bigg(&\frac{(X_i+X_j \pm \alpha_n(i,j))^2}{4\beta}\bigg)\\
& =\exp\bigg(\frac{(X_i+X_j)^2}{4\beta}\bigg)\exp\bigg(\frac{\alpha^2_n(i,j)}{4\beta}\bigg)\exp\bigg( \pm \frac{2(X_i+X_j)\alpha_n(i,j)}{4\beta}\bigg)\\
& \leq \exp\bigg(\frac{(X_i+X_j)^2}{4\beta}\bigg)\exp\bigg(\frac{C^2_n}{4\beta}\bigg)\exp\bigg( \frac{2D_n}{4\beta}\bigg).
\end{split}
\end{align*}
From extreme value theory (see, e.g. \cite{galambos}, p. 227) we have
$\max_{1\leq i\leq n} |X_i| = O_\PP\left(\sqrt{\log(n)}\right)$. Since $C_n$ and $D_n$ do not depend on $i$ and $j$, it follows that
\[ C_n =O_\PP\left(\frac{\sqrt{\log(n)}}{\sqrt{n}}\right) = o_\PP(1), \quad  D_n =O_\PP\left(\frac{\log(n)}{\sqrt{n}}\right) = o_\PP(1)
\]
and thus
\[
\exp\bigg(\frac{C^2_n}{4\beta}\bigg)\exp\bigg( \frac{2D_n}{4\beta}\bigg) = 1 + o_\PP(1).
\]
Consequently,
\begin{align*}
\|R_{n,1}\|^2 &\leq O_\PP\left(n^{-1}\right) \frac{\sqrt{\pi}}{2\beta^{9/2}}\frac{1}{n^2}\sum_{i,j =1}^{n}\Bigg(X_i^2X_j^2\\
& \hspace{3.5cm}\cdot\big( (X_i+X_j)^4 + \alpha_n^4+ 6\beta(X_i+X_j)^2 + 6\alpha_n^2 + 3\beta^2\big)\\
&\hspace{3.5cm} \cdot \exp\bigg(\frac{(X_i+X_j)^2}{4\beta}\bigg)\Bigg)\exp\bigg(\frac{C^2_n}{4\beta}\bigg)\exp\bigg( \frac{2D_n}{4\beta}\bigg)\\
& = O_\PP\left(n^{-1}\right)\frac{\sqrt{\pi}}{2\beta^{9/2}}\frac{1}{n^2}\sum_{i,j =1}^{n}\Bigg(X_i^2X_j^2\big( (X_i+X_j)^4 + 6\beta(X_i+X_j)^2 + 3\beta^2\big)\\
&\hspace{5cm}\cdot\exp\bigg(\frac{(X_i+X_j )^2}{4\beta}\bigg)\Bigg)(1+o_\PP(1)),
\end{align*}
and
\begin{align*}
\|R_{n,2}\|^2 &\leq O_\PP\left(n^{-1}\right) \frac{\sqrt{\pi}}{2\beta^{9/2}}\frac{1}{n^2}\sum_{i,j =1}^{n}\Bigg(\big( (X_i+X_j)^4 + \alpha_n^4+ 6\beta(X_i+X_j)^2 + 6\alpha_n^2 + 3\beta^2\big)\\
&\hspace{3cm} \cdot \exp\bigg(\frac{(X_i+X_j)^2}{4\beta}\bigg)\Bigg)\exp\bigg(\frac{C^2_n}{4\beta}\bigg)\exp\bigg( \frac{2D_n}{4\beta}\bigg)\\
& = O_\PP\left(n^{-1}\right)\frac{\sqrt{\pi}}{2\beta^{9/2}}\frac{1}{n^2}\sum_{i,j =1}^{n}\Bigg(\big( (X_i+X_j)^4 + 6\beta(X_i+X_j)^2 + 3\beta^2\big)\\
&\hspace{5cm}\cdot\exp\bigg(\frac{(X_i+X_j )^2}{4\beta}\bigg)\Bigg)(1+o_\PP(1)).
\end{align*}
Since $\beta >2$ we have
\begin{align*}
\BE\bigg[X_1^2X_2^2\big((X_1+X_2)^4+6\beta(X_1+X_2)^2 + 3\beta^2\big)\exp\bigg(\frac{(X_1+X_2)^2}{4\beta}\bigg)\bigg] &< \infty,\\
\BE\bigg[X_1^4\big(16X_1^4+24\beta X_1^2 + 3\beta^2\big)\exp\bigg(\frac{X_1^2}{\beta}\bigg)\bigg]&< \infty,\\
\BE\bigg[\big((X_1+X_2)^4+6\beta(X_1+X_2)^2 + 3\beta^2\big)\exp\bigg(\frac{(X_1+X_2)^2}{4\beta}\bigg)\bigg] &< \infty,\\
\BE\bigg[\big(16X_1^4+24\beta X_1^2 + 3\beta^2\big)\exp\bigg(\frac{X_1^2}{\beta}\bigg)\bigg]& < \infty,
\end{align*}
and hence
\begin{align*}
\frac{1}{n^2}\sum_{i,j =1}^{n}\left(X_i^2X_j^2\big( (X_i+X_j)^4 + 6\beta(X_i+X_j)^2 + 3\beta^2\big)\exp\bigg(\frac{(X_i+X_j )^2}{4\beta}\bigg)\right) &= O_\PP(1),\\
\frac{1}{n^2}\sum_{i,j =1}^{n}\left(\big( (X_i+X_j)^4 + 6\beta(X_i+X_j)^2 + 3\beta^2\big)\exp\bigg(\frac{(X_i+X_j )^2}{4\beta}\bigg)\right) &= O_\PP(1).
\end{align*}
Summarizing, it follows that $\|R_{n,1}\|^2 \leq O_\PP\left(n^{-1}\right)$ and $\|R_{n,2}\|^2 \leq O_\PP\left(n^{-1}\right)$ and thus
$\|R_n\|^2 = o_\PP(1)$, which is (\ref{asymneg}).

\section{Remarks and open problems}
\subsection{Remark (An alternative approach via V-statistics)} \vspace*{-3mm}
Under more restrictive conditions on $\beta$, the limit null distribution of $T_{n,\beta}$ may also be obtained using results of
\cite{dewet}. To this end, let $\vartheta =(\mu, \sigma^2) \in \Theta := \RR \times \RR_{>0}$ and put
\begin{align*}
h_\beta(x,y;\vartheta)& := \sqrt{\pi}\Bigg(\frac{1}{\sqrt{\beta-1}} -\frac{1}{\sqrt{\beta-\frac{1}{2}}}\bigg( \exp\bigg(\frac{(x-\mu)^2}{(4\beta-2)\sigma^2}\bigg) + \exp\bigg(\frac{(y-\mu)^2}{(4\beta-2)\sigma^2}\bigg) \bigg) \\
&+ \frac{1}{\sqrt{\beta}} \exp\bigg(\frac{(x+y-2\mu)^2}{4\beta\sigma^2}\bigg)\Bigg).
\end{align*}
Letting $\widehat{\vartheta}_n=(\overline{X}_n,S_n^2)$, we have
\[
\frac{T_n}{n} = \frac{1}{n^2}\sum_{i,j=1}^{n}h_\beta(X_i,X_j;\widehat{\vartheta}_n),
\]
which means that $T_{n,\beta}/n$ is a V-statistic with estimated parameters. Moreover, putting
\[
g(x,t;\vartheta) := \exp\bigg(\frac{t(x-\mu)}{\sigma}\bigg)-\exp\bigg(\frac{t^2}{2}\bigg), \quad x,t \in \RR,
\]
we have
\[
h_\beta(x,y;\vartheta) = \int_\RR g(x,t;\vartheta)g(y,t;\vartheta)\exp(-\beta t^2) \, \dd t,
\]
which shows that $T_{n,\beta}/n$ is the special type of $V$-statistic considered in \cite{dewet}.

\subsection{Remark (Contiguous alternatives)} \vspace*{-3mm}
Suppose $X_{n,1},\ldots,X_{n,n}$ are i.i.d. random variables with the density
\begin{equation}\label{contig}
f_n(x) = \varphi(x) \left( 1 + \frac{g(x)}{\sqrt{n}} \right),
\end{equation}
where $\varphi$ is the density of the standard normal distribution and $g: \RR  \rightarrow \RR$ is a bounded measurable function
satisfying $\int_\RR g(x) \varphi(x) \, \dd x = 0$. We assume that $n$ is suffiently large to ensure that $f_n$ is nonnegative. Put
\[
c(t) := \int_\RR h(x,t) g(x) \varphi(x) \, \dd x, \quad t \in \RR,
\]
where $h(x,t)$ is given in (\ref{hvonxt}), and let $P_n := \otimes_{j=1}^n (\varphi \lambda^1)$, $Q_n := \otimes_{j=1}^n (f_n \lambda^1)$,
where $\otimes$ denotes product measure and $\lambda^1$ is Borel Lebesgue measure on ${\cal B}$. Putting $L_n = \dd Q_n/ \dd P_n$,
we have
\[
\log L_n = \sum_{j=1}^n \log \left( 1+ \frac{g(X_{n,j})}{\sqrt{n}} \right) = \sum_{j=1}^n \left( \frac{g(X_{n,j})}{\sqrt{n}} - \frac{g^2(X_{n,j})}{2n}\right) + o_{P_n}(1)
\]
and thus, by the Central Limit Theorem and Slutzki's Lemma
\[
\log L_n \, \vertk \, {\rm N}\left( - \frac{\sigma^2}{2}, \sigma^2 \right) \quad {\rm{ under }} \ P_n,
\]
where $\sigma^2 = \int_\RR g^2(x) \varphi(x) \, \dd x$.
Invoking LeCam's first lemma (see, e.g., \cite{witting}, p. 311), the sequence $Q_n$ is contiguous to $P_n$.
Straightforward algebra shows that, under $P_n$,
\[
\lim_{n \to \infty} \Cov(\widetilde{W}_n(t), \log L_n) =c(t),
\]
where $\widetilde{W}_n$ is the process defined in (\ref{wntilde}).
Therefore, for fixed $k$ and $t_1,\ldots,t_k \in \RR$,
the joint limiting distribution of
$\widetilde{W}_n(t_1), \ldots , \widetilde{W}_n(t_k)$ and $\log L_n$ under $P_n$, as $n \to \infty$,  is the ($k+1$)-variate normal distribution
\[
{\rm N}_{k+1} \left(
\begin{pmatrix} 0 \\ \vdots \\ 0 \\ - \frac{\sigma^2}{2} \end{pmatrix},
\begin{pmatrix} \Sigma & \mathbf{c} \\
\mathbf{c}^\top & \sigma^2
\end{pmatrix}
\right),
\]
where $\Sigma = (K(t_i,t_j))_{1 \le i,j \le k})$ with $K$ given in (\ref{kernel}) and $\mathbf{c} = (c(t_1), \ldots, c(t_k))^\top$.
By LeCam's third lemma (see, e.g., \cite{witting}, p. 329), the finite-dimensional
distributions of $\widetilde{W}_n$  converge under $Q_n$ to the finite-dimensional distributions of the
shifted Gaussian element $W+c$, where $W$ is given in Theorem \ref{nulldistri}.
Since tightness of $\widetilde{W}_n$ under $P_n$ and the
contiguity of $Q_n$ to $P_n$ entail tightness of $\widetilde{W}_n$ under $Q_n$, we have
$\widetilde{W}_n \vertk W+c$ under $Q_n$. Since $\|W_n - \widetilde{W}_n\| = o_{P_n}(1)$ (see (\ref{wnwntild})) and thus
 $\|W_n - \widetilde{W}_n\| = o_{Q_n}(1)$  by contiguity, we have $W_n \vertk W+c$ under $Q_n$.
 The Continuous Mapping Theorem then yields
\[
T_{n,\beta} \, \vertk \, \int_{\RR} \left(W(t) + c(t)\right)^2 \, \exp(-\beta t^2) \, \dd t \quad {\rm{ under }} \ Q_n \ {\rm{ as }} \ n \to \infty.
\]
Thus, $T_{n,\beta}$ has a limit distribution under contiguous alternatives to $H_0$ given by (\ref{contig}).

\subsection{Remark (Two open problems)}
Denoting the right-hand side of (\ref{unglie}) by $\Delta$, we conjecture that
\[
\frac{T_{n,\beta}}{n} \, \rightarrow \, \Delta \quad \textrm{in probability as } n \to \infty.
\]
Such a result would open the ground for tackling asymptotic normality of
\[
\sqrt{n}\left(\frac{T_{n,\beta}}{n}  - \Delta \right)
\]
under fixed alternatives as $n \to \infty$, in the spirit of \cite{bahe}.

Regarding consistency, we conjecture that $\lim_{n\to \infty} T_{n,\beta} = \infty$
$\PP$-almost surely under {\em any} fixed alternative distribution. Hence, the test based on $T_{n,\beta}$ would be
globally consistent.

\bigskip


N. Henze, Institute of Stochastics, Karlsruhe Institute of Technology (KIT), Englerstr. 2, D-76133 Karlsruhe:
\\
{\texttt Norbert.Henze@kit.edu}
\vspace*{2mm}


S. Koch, Institute of Mathematics, University of Mannheim, A5 6,  D-68159 Mannheim:
\\
{\texttt stefan.koch@uni-mannheim.de}

\end{document}